\documentclass[10pt]{article}
\usepackage{graphicx}

\def\R{{\rm I\! R}}
\def\N{{\rm I\! N}}

\def \pTW*{\partial_{W^*} T}
\def \pTnW*{\partial^{(n)}_{W^*} T}
\def \ps{\frac{\partial }{\partial s}}
\def \ds{\frac{d }{d s}}
\def \dns{\frac{d^{n} }{d s^n}}
\def \pns{\frac{\partial^{n} }{\partial s^n}}

\newtheorem{theorem}{Theorem}
\newtheorem{corollary}{Corollary}

\newtheorem{lemma}{Lemma}

\newtheorem{remark}{Remark}

\title{The period functions' higher order derivatives}

\author{M. Sabatini
\footnote{Dip. di Matematica, Univ. di Trento, I-38050 Povo, (TN) - Italy.
Phone: ++39(0461)281670, Fax: ++39(0461)281624, Email: marco.sabatini@unitn.it - \ \ \ \ \ \ \ \  \ \ \ \ \ \ \ \ \ \ \ \ This paper was partially supported by the PRIN project {\it Equazioni differenziali ordinarie: sistemi dinamici, metodi topologici e applicazioni}. }
}
\date{February 9$^{th}$, 2012}
\begin{document}
\maketitle
\begin{abstract} We prove a formula for the $n$-th derivative of the period function $T$ in a  period annulus of  a planar differential system. For  $n = 1$, we obtain Freire, Gasull and Guillamon formula for the period's first derivative \cite{FGG}. We apply such a result to hamiltonian systems with separable variables and other systems. We give some sufficient conditions for  the period function of conservative second order O.D.E.'s  to be convex. 

\end{abstract}

\section{Introduction}

Let $\Omega$ be an open connected subset of the real plane. Let us consider a differential system
\begin{equation}\label{sysV} 
z' = V(z), \qquad z\equiv (x,y) \in \Omega,
\end{equation} 
$V(z)=  (V_1(z),V_2(z))\in C^\infty(\Omega,\R^2)$. We denote by $\phi_V(t,z)$ the local flow defined by (\ref{sysV}). A topological annulus $A \subset \Omega$ is said to be a {\it period annulus} of (\ref{sysV}) if it is the set-theoretical union of concentric non-trivial cycles of (\ref{sysV}). If the inner component of $A$'s boundary  is a single point $O$, then $O$ is said to be a {\it center}, and the largest connected punctured neighbourhood $N_O$ of $O$ covered with non-trivial cycles is said to be its {\it central region}. If $A$ is a period annulus, we can define on $A$ the {\it period function T} by assigning to each point $z\in A$ the minimum positive period $T(z)$ of the cycle $\gamma$ passing through $z$. We say that the period function $T$ is {\it increasing} if outer cycles have larger periods.  Let $\gamma^\dashv(s)$ be a curve of class $C^1$ meeting transversally the cycle $\gamma$ at the point $s=s_0$.
We say that $\gamma$ is a {\it critical cycle} if $\left[ \frac{d}{ds}T(\gamma^\dashv(s))\right]_{s=s_0} =0 $. It is possible to prove that such a definition does not depend on the particular transversal curve $\gamma^\dashv$ chosen.  

The existence and number of critical cycles affects the number of solutions to some boundary value problems. In fact, given a positive $\tau \in \R$, the number of $\tau$-periodic cycles contained in a period annulus $A$ is bounded above by $n+1$, where $n$ the number of critical cycles  contained in $A$. Similarly, the existence of critical orbits is related to the study of some Neumann problems for systems equivalent to second order differential equations, as well as to the study of mixed problems. The absence of critical orbits is itself an important element in the treatment of  boundary value problems, bifurcation or perturbation problems \cite{ChJ, Sch}, delay differential equations  \cite{CGM}, thermodynamics \cite{R2}, linearizability \cite{MRT}. 

The simplest case, that of monotone period functions, was dealt with in several papers. For several references we address to the bibliographies of  \cite{FGG,  Vt, YZ, Z}. A very special sub-case is that of isochronous systems, i. e. systems having period annuli with constant period function, for which we refer to the survey \cite{CS}, and to the bibliographies of some recent papers as  \cite{BCS, LV}. Systems with one or more critical orbits were studied in \cite{BBLT, CLH, ChD1, CGM, GGJ1, GGJ2, G, R2, S3, W, YZ, Z}. 

In some papers upper bounds to the local or global number of critical orbits are studied \cite{ChD2, ChJ, ChS}. Finding an upper bound to the number of critical cycles is similar, to some extent, to the problem of finding an upper bound to the number of limit cycles of a planar differential system. In fact, similar techniques have been developped in the treatment of such problems, mainly in a bifurcation perspective. In particular, the role played by the displacement function's derivatives in limit cycles' bifurcation is similar to that of the period function's derivatives for critical cycles bifurcation  \cite {ChD2, ChJ, CGdS, F96,F98, MMV}. Dealing with bifurcation from a critical point $O$, the key piece of information is the order of the first non-vanishing derivative at $O$. On the other hand, global estimates to the number of critical orbits are usually based on some property of $T$'s derivatives in all of a period annulus. This task has been faced in different ways in \cite{BBLT, CLH, ChS, G, LL, MV, W} for second order ODE's,  \cite{ChJ, CG, S3} for other types of hamiltonian systems,  \cite{GGJ1, GGJ2} for complex differential equations, \cite{GLY} for polynomial systems. 

In this paper's view, a key result  was given in \cite{FGG}. Let us
denote by $[V,W] = \partial_V W -\partial_W V$ the Lie bracket of $V$ and $W$. A vector field $W$, transversal to $V$, is said to be a {\it non-trivial normalizer} of $V$ on a set $A\subset \Omega$ if there exists a function $\mu$ defined on $A$ such that $[V,W] = \mu V$ on $A$. We call $\mu$ its {\it N-cofactor}. 
 Let $\phi_W(s,z)$ be the local flow defined by the solutions of
\begin{equation}\label{sysW} 
z' = W(z) .
\end{equation} 
In \cite{FGG},  it was proved that 
\begin{equation}\label{FGG}
\partial_W T(z) = \frac{d}{ds} \ T(\phi_W(s,z)) = \int_0^T \mu(\phi_V(t,z)) dt.
\end{equation} 
Hence, if a normalizer is known, such an approach allows to get some information about the first derivative of the period function, in particular when the function $\mu$ does not change sign, avoiding the need to evaluate the above integral. In order to find a normalizer, it is sufficient to know a first integral, as shown in  \cite{S3}. In some special cases, more convenient normalizers can be found, as for Hamiltonian systems with separable variables \cite{FGG}.

In this paper we give a formula for the $n$-th derivative of the period function, based on an approach similar to that one introduced in  \cite{FGG}. Rather than focusing on a $N$-cofactor, we base our approach on a the existence of a suitable commutator, looking for a function $m(z)$ such that $[mV,W] = 0$. We call $m$ a {\it C-factor}. C-factors and N-cofactors are related by the equality $\partial_W m = m \mu$, so that they can be obtained from each other, at least in principle.
Let us denote by $\partial_W^{(n)} T$ the $n$-th derivative of $T(\phi_W(s,z))$ with respect to $s$. In our main result, we prove that
\begin{equation}
\partial_W^{(n)} T(z) = \frac{d^n}{ds^n} \ T(\phi_W(s,z)) =  \int_0^{T(z)} \frac {\partial_W^{(n)} m(\phi_V(t,z)) }{m(\phi_V(t,z))} dt,
\end{equation}
In the case of the first derivative, our formula reduces to  (\ref{FGG}). We also give a recursive formula that, starting from a N-cofactor $\mu$, allows to avoid C-factors. Setting
$$
\mu_1 = \mu, \qquad \mu_n = \mu_{n-1} \mu  + \partial_W \mu_{n-1},
$$
we prove that
$$
\partial_W^{(n)} T(z) = \int_0^{T(z)} \mu_n(\phi_V(t,z))   dt .
$$
In particular, $T$ is $W$-convex if $\mu_2 = \mu^2  + \partial_W \mu \geq 0$. 
We also provide a wide class of systems having explicit C-factors, including hamiltonian systems with separable variables. In such a case, we extend the formula given in \cite{FGG} for hamiltonian systems with separable variables to higher order derivatives. In particular, we prove that a hamiltonian system with separable variables
$$
x' = F'(y), \qquad  y' = - G'(x) 
$$ 
has at most one critical orbit if
$$
\mu_{s2} = 4 \bigg[  1   +  2\,\frac{ G  G'' }{G'^2 } \cdot\frac{ F  F''}{F'^2 } +
$$ $$
+  \frac { 3 G^2 G''^2 - 3 G G'^2 G''  -  G^2 G' G''' }{G'^4} + \frac{ 3  F^2 F''^2 - 3 F F'^2 F'' -   F^2 F' F''' }{F'^4} \bigg].
$$
has constant sign. As a consequence, we prove that for systems equivalent to second order differential equations, the inequality 
$$
 G'^4 - 8 G G'^2 G''  +12 G^2 G''^2 - 4 G^2 G' G'''   \geq 0
$$
implies $T$'s convexity.

\section{Some properties of normalizers and period functions}

We assume $V$ and $W$ to be transversal on  a period annulus $A$, and $W$ to be a non-trivial normalizer of $V$, i. e. $V \wedge W \neq 0$, and $[V,W]=\mu V$ on $\Omega$. One has 
$\displaystyle{  \mu = \frac{[V,W] \cdot V}{|V|^2} \in C^\infty(A,\R)}$. 
Let us  choose arbitrarily a point $z^*$ in $A$. 
The orbit $\phi_W(s ,z^*)$ meets all the cycles of $A$. Since $W$ is  a normalizer, the  map $\phi_W(s ,\cdot)$ takes $V$-cycles into $V$-cycles, hence there is a one-to-one correspondence between $V$-cycles and the values of $s$ in some real interval $I_W$. 
We parametrize such cycles by means of the parameter $s$. 
The transversality of $V$ and $W$ implies that $T\in C^\infty (A,\R)$.  Since, by definition, $T$ is constant on cycles, hence a first integral of (\ref{sysV}), we may write $T(s)$ to denote the period of the unique $V$-cycle corresponding to the parameter $s$. Different normalizer can produce the same parametrization. In fact, since, by the transversality of $V$ and $W$ every vector field definied in $A$ can be expressed as a linear combination $Z = \alpha V + \beta W$, $\beta >0$ in order to preserve transversality, one has
$$
[V,Z] = [V,\alpha V + \beta W] = (\beta \mu + \partial_V \alpha ) V +  (\partial_V \beta) W .
$$
Hence $\alpha V + \beta W$ normalizes $V$ if and only if  $\partial_V \beta = 0$, i. e. $\beta$ is a first integral of $(\ref{sysV})$. Moreover, the new N-cofactor is $\beta \mu + \partial_V \alpha$. In a sense, we can split the action of $\alpha V + \beta W$ on $A$ as the combination of a rotation along $V$-orbits, determined by the term $\alpha V$, plus a motion along the $W$-orbits determined by the term $\beta W$. Applying the main theorem in \cite{FGG} that  gives the first derivative of $T$ w. resp. to $s$, 
$$
\partial_Z T  = \int_0^T (\beta \mu + \partial_V \alpha)(\phi_V(t,z))\ dt = 
$$
$$ 
= \beta \int_0^T \mu (\phi_V(t,z))\ dt = \beta\ \partial_W T,
$$
we see that the contribution of the term $\partial_V \alpha$ is zero, while the presence of the $V$-first integral $\beta$ only appears as a factor that can be taken out of the integral in (\ref{FGG}), since it is constant with respect to $t$ on $\phi_V(t,z)$.
In conclusion, every normalizer gives the same expression for $T'(s)$, up to a multiplicative factor, which is a first integral of $(\ref{sysV})$.

In some cases a particular parametrization for the $V$-cycles is preferred. This is the case of hamiltonian systems
\begin{equation}\label{sysham}
x' = H_y, \qquad y' = -H_x,
\end{equation}
where the \lq\lq natural" parameter is provided by the Hamiltonian function $H(z)$. A normalizer producing $H(z)$ as a cycle's parameter is
\begin{equation}\label{syshamnorm}
x' = \frac{H_x}{H_x^2 + H_y^2}, \qquad y' = \frac{H_y}{H_x^2 + H_y^2},
\end{equation}
since for such a system one has $\dot H = 1$. In \cite{S4} it was proved that the related N-cofactor is 
\begin{equation}\label{muham}
\mu_H =   
\frac{( H_{yy} - H_{xx})H_x^2 - 4H_{xy}H_xH_y  + (H_{xx}  - H_{yy})H_y^2}{|\nabla H|^4} ,
\end{equation}
so that 
\begin{equation}\label{TpH}
T'(H) =  \int_0^{T(H)} \mu_H(\phi_V(t,z)) dt ,
\end{equation}
performing the integration along a  $V$-cycle $\phi_V(t,z)$. For some hamiltonian systems it is possible to find other normalizers such that $\dot H = 1$. For instance, if $H(x,y) = F(y) + G(x)$, one can consider the system
$$
x' = \frac{G}{H G'}, \qquad y' = \frac{F}{H F'},
$$
a special case of  the normalizer $U^*$ in \cite{FGG}, remark 7, obtained choosing $\displaystyle{  k(x,y) =   }$ $\displaystyle{  = \frac{1}{H(x,y) }  }$. On the other hand, such a normalizer does not exist at points $z$ such that $G(z) \neq 0$, $G'(z) = 0$. 

If $W$ is a non-trivial normalizer of $V$,  then it is a non-trivial normalizer of $\rho V$, for every non-vanishing function $\rho \in C^\infty(\Omega,\R^2)$. In fact, assuming $\rho > 0$, one has
$$
[\rho V, W] = (\rho \mu - \partial_W \rho) V = \Big (\mu  - \partial_W (\ln \rho) \Big) (\rho V) = \overline \mu \, (\rho V). 
$$
This finds an application when a first integral $H(x,y)$ of a system (\ref{sysV}) is known, since the system (\ref{syshamnorm}) is a normalizer of every system having $H(x,y)$ as a first integral.  We say that a non-vanishing function $\rho \in C^\infty(\Omega,\R^2)$ is an {\it inverse integrating factor} of the system (\ref{sysV}) if  $(-V_2(z),V_1(z) )= \rho(z) \nabla H(z)$, for some $H(x,y)$. If this occurs, then the N-cofactor corresponding to the normalizer (\ref{syshamnorm}) is 
$$
\overline \mu_H = \mu_H - \partial_W (\ln \rho). 
$$
\medbreak 
Next lemma shows that differentiating with respect to $s$ produces new first integrals of (\ref{sysV}). Let us set

\begin{lemma} For every $n \geq 1$, $T^{(n)}$ is a first integral  of $(\ref{sysV})$.
\end{lemma}
{\it Proof.} it is sufficient to observe that both $s$  and $T(s)$ are first integrals of (\ref{sysV}), so that every derivative of the period function w. resp. to $s$ depends only on $s$, i.e. it is constant on the cycles of (\ref{sysV}).  \hfill$\clubsuit$\medbreak

In the following we shall denote by
$$
T^{(n)}(z) = \partial_W^{(n)} T = \frac {d^{n} }{ds^n} T(\phi_W(s ,z))
$$
the $n$-th derivative of $T$ with respect to $s$. 
Dealing with $T^{(n)}$, for $n > 1$, leads to change one' approach in relationship to the choice of normalizers. In fact, when studying the sign of $T'$ it makes no difference to use any normalizer, but for computational complexity. On the other hand, different normalizers might produce higher order derivatives with different signs.  Let us denote by $\sigma$ the parameter induced by $Z$, so that $\phi_W(s,z) = \phi_Z(\sigma(s),z)$, with $\sigma'(s) >0$. One has
$$
\frac {d }{ds} T(\phi_W(s ,z)) =\sigma'(s) \frac {d }{ds} T(\phi_Z(\sigma(s) ,z)) ,
$$
while
$$
\frac {d^2 }{ds^2} T(\phi_W(s ,z)) =\sigma''(s) \frac {d }{ds} T(\phi_Z(\sigma(s) ,z)) + ( \sigma'(s) )^2\frac {d^2 }{ds^2} T(\phi_Z(\sigma(s) ,z)) .
$$
This shows that convexity is not normalizer-invariant, so that we shall say that a period function $T$ is {\it $W$-convex}, rather  than just convex. This opens the problem to find the most convenient normalizer $W$ in order to study the sign of $\partial_W^{(n)} T$. \medbreak

In general, the unique  $\partial_W^{(n)} T$ to have normalizer-invariant sign is $\partial_W T$. If a cycle is critical with respect to a normalizer, then it is critical with respect to any other normalizer, and a period annulus $A$ can be split into sub-annuli where $T$ is monotonic with respect to any normalizer. 

\medbreak

In this paper's applications, we mainly consider the convexity problem and its consequence on critical cycle's uniqueness, in particular around a center $O$. This leads to different situations for degenerate centers and for non-degenerate ones. 
In order to illustrate the behaviour of the period function in a neighbourhood of a center, let us restrict to analytic systems. If $O$ is a non-degenerate center, then $T$ has an analytic extension at the origin \cite{V}. 
Since analytic functions cannot have infinitely many zeroes accumulating at a point, every line segment  passing through the origin has an open subsegment  $\Sigma$ such that all cycles meeting $\Sigma$ have only two points on $\Sigma$. Possibly rotating the axes, we may assume $\Sigma$ to be contained in the $x$-axis, i.e. $\Sigma =\{ (x,0): -\varepsilon < x < \varepsilon,  \varepsilon > 0 \}$. Then one can define an involution $\iota: \Sigma \rightarrow \R$, i. e. a function such that  such that $\iota( \iota(x))=x$, satisfying $T(x) = T(\iota(x))$ \cite{MV}. For $x$-reversible centers, such a property reduces to $T(x) = T(-x)$, i. e. $T$ is even. 
As a consequence, its Taylor expansion on the $x$-axis has the form
$$
T(x,0) = a_{2k} x^{2k} + a_{2k+1} x^{2k+1} + ...
$$
Hence, if $a_{2k} >0$, it is both increasing and convex (decreasing and concave) in a neighbourhood of $O$.  Hence, if we want to prove the uniqueness of a critical orbit, we can only try to prove convexity/concavity out of a neighbourhood of $O$. 
If the center is degenerate, then the period function is unbounded at the origin, so that $T$ could be both decreasing in a neighbourhood of $O$ and convex.

\section{Results}

We say that a map {\it linearizes} a differential system if it takes such a system into a linear one. Next lemma was contained in the unpublished preprint \cite{Sp} as theorem 3.

\begin{lemma} \label{lemmalin}
Let $A$ be a period annulus of (\ref{sysV}), with $[V,W] =0$ on $A$. Then there exists a map $\Lambda \in C^\infty(A,\R^2)$ that linearizes both (\ref{sysV}) and (\ref{sysW}). 
\end{lemma}
{\it Proof.} 
By hypothesis, $V$ and $W$ commute, i. e.  if $\phi_V(t, (\phi_W(s,z)))  $ and  $ \phi_W(s, (\phi_V(t,z)))  $ exist for all $s$ and $t$ in a rectangle $I_s \times I_t$, , then $\phi_V(t, (\phi_W(s,z)))$ $  = \phi_W(s, (\phi_V(t,z)))  $. 
By the main theorem in \cite{S0}, $A$ is an isochronous annulus, i. e. every $V$-cycle in $A$ has the same period $ T>0$. 
Possibly multiplying $V$ by $\frac{T}{2\pi}$, we may assume (\ref{sysV}) to have period $2\pi$.  Following  \cite{S0}, we choose arbitrarily a point $z^*$ in $A$ and define on $A$ the functions $t(z)$, $s(z)$ such that $\phi_V(t(z), (\phi_W(s(z),z^*))) = \phi_W(s(z), (\phi_V(t(z),z^*)))  = z$. The regularity of $t(z)$, $s(z)$ can be proved by the implicit function theorem, as in  \cite{S0} or \cite{FSZ}, section 4. 
By construction, one has 
\begin{equation} \label{tsVW}
 \left\{ \begin{array}{rl}
\partial_V s(z) &= 0 \\
\partial_V t(z) &= 1,
\end{array} \right.
\qquad\qquad\qquad
 \left\{\begin{array}{rl}
\partial_W s(z) &= 1 \\
\partial_W t(z) &=   0.
\end{array} \right.
\end{equation} 
Let us define $\Lambda: A \mapsto \R^2 $ as follows:
$$
 (u,v) = \Lambda (z) = (e^{s(z)}\sin t(z),e^{s(z)}\cos t(z)) .
$$
In order to show that $\Lambda$ is injective, first consider that there is a one-to-one correspondence between $V$-cycles and values of $s(z)$. Moreover, there is a one-to-one correspondence between points of a $V$-cycle and the values of $t(z)$, for $t(z) \in [0,2\pi)$. The map $\Lambda$ transforms injectively the $V$-cycle corresponding to $s(z)$ into the circle $(e^{s(z)}\sin t(z),e^{s(z)}\cos t(z))$.

Using (\ref{tsVW}) it is immediate to show that $\Lambda$ transforms (\ref{sysV}) and (\ref{sysW}) respectively into the following systems
\begin{equation} \label{sysuv}
 \left\{\begin{array}{rl}
u' &= v \\
v' &= -u
\end{array} \right. ,
 \qquad\qquad\qquad
\left\{\begin{array}{rl}
u' = u \\
v' = v
\end{array} \right. .
\end{equation} 
\hfill$\clubsuit$\medbreak

In next theorem we give a formula for the $n$-th derivative of $T$ w. resp. to the parametrization induced by a transversal vector field $W$ such that $[mV,W]= 0$, for some non-vanishing multiplier $m$. We may assume $m$ to be positive, since the period function of $V$ and that one of $-V$ coincide.

\begin{theorem} \label{teordern} Let $A$ be a period annulus of (\ref{sysV}) and $m\in C^\infty(A,\R)$, $m \neq 0$,  such that $[mV,W] =0 $ on $A$. Then
\begin{equation}  \label{dern}
\partial_W^{(n)} T(z) = \int_0^{T(z)} \frac {\partial_W^{(n)} m(\phi_V(t,z)) }{m(\phi_V(t,z))} dt,
\end{equation}
\end{theorem}
{\it Proof.} 
by hypothesis, the systems
$$
z' = m(z)V(z) \qquad z' = W(z) 
$$
commute. By lemma \ref{lemmalin}, there exists a transformation $\Lambda\in C^\infty(A,\R^2)$ which takes such systems resp. into the systems (\ref{sysuv}). As a consequence, the system (\ref{sysV}) is transformed into the system
\begin{equation} \label{V*}
 \left\{\begin{array}{rl}
u' &= \frac{v}{m^*(u,v)} \\
v' &= -\frac{u}{m^*(u,v)} 
\end{array} \right. ,
\end{equation}
where $m^*(u,v) = m(\Lambda^{-1}(u,v))$. Since $\Lambda$ preserves the time, the system (\ref{V*}) has a period annulus with cycles having the same periods as their anti-images in $A$. Denoting by $\theta^*$ the argument function in the plane $(u,v)$, one has
\begin{equation} \label{T}
T = \int_0^T 1\ dt = \int_0^{2\pi} \frac {d\theta^*}{\dot \theta^*} = - \int_0^{2\pi} m^* d\theta^* .
\end{equation} 
The system (\ref{sysW}) is transformed into the second system in (\ref{sysuv}), which is a normalizer of (\ref{V*}). We denote its vector field by $W^*$. We can find the derivative of $T$ with respect to the parameter $s$ by differentiating the integral in (\ref{T}) with respect to $s$:
\begin{equation} \label{Tenne}
\pTW* = \ds T   = - \ds \int_0^{2\pi} m^* d\theta^* = - \int_0^{2\pi} \ps m^* d\theta^* =  \int_0^T \frac{1}{m^*} \ps m^*\ dt.
\end{equation} 
As for higher order derivatives, a similar conclusion holds, since the integration extremes do not depend on $s$:
\begin{equation} \label{Tn}
\pTnW* = \dns T   = - \dns \int_0^{2\pi} m^* d\theta^* = - \int_0^{2\pi} \pns m^* d\theta^* =  \int_0^T \frac{1}{m^*} \pns m^*\ dt.
\end{equation} 
Since $\Lambda$ preserves the time of all the systems considered, applying the inverse transformation $\Lambda^{-1}$ gives the formula 
(\ref{dern}).
\hfill$\clubsuit$\medbreak

In next lemma we show that the relationship $[mV,W]=0$, for a non-vanishing $m$, is equivalent to $W$ being a normalizer of $V$, and find the relationship between $m$ and the N-cofactor $\mu$. We state it for a period annulus, even if it holds in other subsets of $\Omega$.

\begin{lemma} \label{lemmam} Let $V, W \in C^\infty (A,\R^2)$. Then 
\begin{itemize}
\item[i)] there exists $m \in C^\infty(A,\R)$, $m(z) > 0$, such that $[mV,W] =0$, if and only if $W$ is a normalizer of $V$, with cofactor  $\displaystyle{ \mu= \partial_W ( \ln m )}$;
\item[ii)]  $m, \overline m \in C^\infty(A,\R)$, $m(z),\overline m(z) > 0$ in $A$, satisfy $[mV,W] = 0 =[\overline mV,W] $, if and only if there exists $J \in C^\infty(A,\R)$, first integral of (\ref{sysW}), such that $\overline m = J m $. Moreover,
\begin{equation} \label{rappn}
 \frac {\partial_W^{(n)} m}{m} =  \frac {\partial_W^{(n)} \overline m}{\overline m} .
\end{equation}
\end{itemize}
\end{lemma}
{\it Proof.} i) 
If $[mV,W] =0$, then
$$
0 = [mV,W]  = m [V,W] -( \partial_W m ) V, 
$$
hence 
$$
[V,W] = \left(\frac{ \partial_W m }{m}  \right) V = \bigg( \partial_W ( \ln m ) \bigg) V.
$$
Conversely, let us assume there  exists $\mu$ such that $[V,W]= \mu V$. Let us choose arbitrarily a point $\overline z \in A$. Let $t(z), s(z)$ be the functions defined in  in lemma  \ref{lemmalin}, so that $z = \phi_V(t(z),\phi_W(s(z),\overline z)) = \phi_W(s(z),\phi_V(t(z),\overline z))$. Let us define $m(z)$,  as follows:
$$
m(z) =e^{ \int_0^{s(z)} \mu(\phi_W(\sigma,\phi_V(t(z),\overline z)) d\sigma    }.
$$
One has $m\in C^\infty(A,\R)$ and
$$
\partial_W m(z) = \mu(z)\, e^{ \int_0^{s(z)} \mu(\phi_W(\sigma,\phi_V(t(z),\overline z)) d\sigma    } = \mu(z) m(z),
$$
since $\partial_W s(z) = 1$. Then
$$
[mV,W] = m [V,W] - (\partial_W m ) V = m \mu V - \mu  m V = 0,
$$

ii) Let $J$ be a first integral of  (\ref{sysW}). Then 
$$
[JmV,W] = Jm [V,W] - (\partial_W (Jm) ) V  = J\Big( m \mu  - \partial_W m  \Big) V  = J(m \mu V - \mu  m V)= 0.
$$
Vice-versa, if $[mV,W] = 0 = [\overline mV,W] $, then one has
$$
 \big( m \mu - \partial_W m  \big)V=   [mV,W]  = 0 = [\overline mV,W] = \big( \overline m \mu - \partial_W \overline  m \big)V.
$$
From  $\mu = \frac{ \partial_W m }{m} =  \frac{ \partial_W \overline m }{\overline m}$ one has
$$
m\ \partial_W  \overline m -  \overline m\  \partial_W  m = 0,
$$
that gives 
$$
\partial_W\left(  \frac{\overline m}{m} \right) =0,
$$
hence $\displaystyle{ \frac{\overline m}{m} }$ is a first integral of (\ref{sysW}).

As for (\ref{rappn}), if $\overline m = J m$, then $\partial_W^{(n)}\ \overline m = \partial_W^{(n)} \big( J m \big) =  J\, \partial_W^{(n)}\,  m $, since $J$ is  a first integral of (\ref{sysW}). Hence the two fractions in (\ref{rappn}) coincide.
\hfill$\clubsuit$\medbreak

In lemma \ref{lemmam} the integration is performed only along the $W$-orbits, starting at the intersection of the $V$-cycle through $\overline z$ with  the $W$-orbit passing through $z$.  Performing the integration along the $W$-orbit, we have chosen an initial value of 1 on the cycle $\phi_V(\cdot,\overline z)$, so that $m(\phi_V(t,\overline z)) = 1$ for all $t\in \R$. We could have chosen other smooth initial values on $\phi_V(\cdot,\overline z)$, getting other isochronous systems with the same period annulus $A$. All such functions $m$ can be obtained from one another multiplying  by a first integral of $W$. All of them give the same ratio that appears in the integral  (\ref{dern}).

As a special case of lemma \ref{lemmam}, we obtain the cited formula for the first derivative of the period function \cite{FGG}.

\begin{corollary} \label{teoFGG} Let $A$ be a period annulus of (\ref{sysV}), with $[V,W] =\mu V$ on $A$. Then
$$
\partial_W T(z) = \int_0^T \mu(\phi_V(t,z)) dt.
$$
\end{corollary}
{\it Proof.} 
It is sufficient to observe that
$$
\partial_W m = \ps m = \ps e^{ \int_0^{s(z)} \mu(\phi_W(\sigma,\phi_V(t(z),\overline z)) d\sigma } = 
\mu(z) e^{ \int_0^{s(z)} \mu(\phi_W(\sigma,\phi_V(t(z),\overline z)) d\sigma}, 
$$
hence
$$
\frac {\partial_W m(z)} {m(z)} = \mu(z). 
$$
\hfill$\clubsuit$\medbreak

In some cases it is computationally more convenient to look for a $\mu$ such that $[V,W] = \mu V$, rather than for a $m$ such that $[mV,W]=0$. In next theorem we give a direct way to compute high-order derivatives of $T$ working only on $\mu$. This can be done by means of a recursive formula.

\begin{theorem} \label{ricorr} 
Let $A$ be a period annulus of (\ref{sysV}), with $[V,W] =\mu V$ on $A$. Then
$$
 \partial^{(n)}_W T(z) = \int_0^T \mu_n(\phi_V(t,z)) dt.
$$
where $\mu_n$ is recursively defined by
\begin{equation}  \label{defricorr}
\mu_1 = \mu, \qquad \mu_n = \mu_{n-1} \mu  + \partial_W \mu_{n-1}.
\end{equation}
\end{theorem}
{\it Proof.} 
Let us define $m$ as in the proof of lemma \ref{lemmam}. Then, let us set
$$
m_n =   \frac{ \partial^{(n)}_W m}{m}.
$$
By corollary  \ref{teoFGG}, one has $m_1 = \mu$.
For $n > 1$, one has 
$$
\partial_W m_{n-1} = 
\partial_W  \left(   \frac{ \partial^{(n-1)}_W m}{m} \right) = 
\frac 1{m^2}  \bigg( m\  \partial^{(n)}_W m -   \big( \partial^{(n-1)}_W  m\big )\ \big( \partial_W m\big) \bigg) =
$$
$$
= \frac {\partial^{(n)}_W m}{m} - \frac {\partial^{(n-1)}_W m}{m}  \frac {\partial_W m}{m} =
m_n - m_{n-1}  m_1,
$$
hence 
$$
m_n = m_{n-1}  m_1 + \partial_W m_{n-1}.
$$
This shows that the functions $m_n$ satisfy the recurrence equations  (\ref{defricorr}), hence $m_n = \mu_n$.
Then the statement comes from theorem  \ref{teordern} .
\hfill$\clubsuit$\medbreak

In the next formulae we replace the notation $\partial_W \mu $ with the simpler one $\mu'$. Similarly for higher order derivatives with respect to $W$. We write the forms of some low order $\mu_n$'s:
$$
\mu_2 = \mu^2 + \mu'.
$$
$$
\mu_3 = \mu^3 + 3 \mu \mu' + \mu''.
$$
$$
\mu_4 = \mu^4 + 6 \mu^2 \mu' + 4 \mu \mu'' + 3 \mu'^2  +\mu'''.
$$
$$
\mu_5 = \mu^5 + 10 \mu^3 \mu' + 10 \mu^2 \mu'' +15 \mu \mu'^2 +5 \mu \mu''' +10 \mu' \mu'' + \mu^{(4)}.
$$

We say that a function $h \in C^\infty(A,\R)$ {\it satisifies the condition  $(B)$} if $h(z)$ does not change sign in $A$, and every cycle contained in $A$ contains a point $z$ such that $h(z) \neq 0$. 

\begin{corollary} \label{cordern} Let $A$ be a period annulus of (\ref{sysV}), with $[mV,W] =0$, with $ m > 0$, and $[V,W] =\mu V$ on $A$. If one of the following holds,
\begin{itemize}
\item[i)]  $\partial_W^{(n)}m$ satisfies the condition $(B)$, 
\item[ii)]  $\mu_n$ satisfies the condition $(B)$, 
\end{itemize}
then $A$ contains at most $n-1$ critical cycles.
\end{corollary}
{\it Proof.}  i) $T$ is a first integral, hence every critical cycle of (\ref{sysV}) corresponds to a critical point of $T(\phi_W(s,\overline z))$. By  (\ref{dern}), if $\partial_W^{(n)}m \geq 0$, then  $\partial_W^{(n)} T \geq 0$. The presence on every cycle of a point where $\partial_W^{(n)} T \neq 0$ implies that on every cycle $\partial_W^{(n)} T > 0$. Similarly if $\partial_W^{(n)}m \leq 0$.

ii) If $\mu_n$ satisfies the condition $(B)$, then also $\partial_W^{(n)}m$ satisfies the condition $(B)$. 
\hfill$\clubsuit$\medbreak

In particular, the corollary  \ref{cordern} holds  for $n=2$, in presence of a convex or concave (w. resp. to the parametrization induced by $W$) period function. The convexity may be proved under a condition on $\partial_W \mu$, rather than on  $ \partial^{(2)}_W m$. 

\begin{corollary}  \label{muconv} Let $A$ be a period annulus of (\ref{sysV}), with $[V,W] =\mu V$ on $A$. If  $\partial_W \mu \geq 0$ on $A$  then $T$ is $W$-convex. If $\partial_W \mu $ satisfies the condition $(B)$, then $T$ is strictly $W$-convex on $A$.
\end{corollary}
{\it Proof.} 
By theorem \ref{ricorr}, one has 
$$
\mu_2 = \mu^2 +  \partial_W \mu \geq 0.
$$
Then
$$
\partial_W^{(2)} T(z) = \int_0^T \mu_2(\phi_V(t,z))   dt \geq 0. 
$$
If every cycle contains a point where $\mu_2 > 0$, the above integral is positive. 
\hfill$\clubsuit$\medbreak

\begin{remark} \label{asym} It is remarkable the asymmetry of such a situation, where convexity can be proved by only considering the sign of $\partial_W \mu$, while concavity cannot. This agrees with the fact that $T$ can be upper unbounded, but it cannot be lower unbounded.
\end{remark}

\section{Applications}

\subsection{Reparametrized isochronous centers}

A center is isochronous if and only if it has a non-trivial commutator \cite{S0}. Given an isochronous center, finding a commutator is not always easy. There are systems which are known to have isochronous centers, for which no commutators are known, as for reversible Li\'enard systems \cite{AFG,CD,CS}. A collection of isochronous systems with their commutators have been provided in \cite{CS}, other isochronous systems may be found in \cite{BCS, LV} and their bibliographies.

In this section we consider centers of the form
\begin{equation} \label{sysxi}
z' = \xi(z) V(z), \qquad \xi \in C^\infty(A,\R), \qquad \xi(z) \neq 0,
\end{equation} 
where $V$ satisfies $[V,W]=0$ on a period annulus $A$, $W$ transversal to $V$, $\xi > 0$. By lemma \ref{lemmam}, an obvious choice to study the period function of $V$  is $\displaystyle{m = \frac 1\xi}$.  Hence we have the following corollary.

\begin{corollary} \label{isorepar} Let $A$ be  a period annulus of ($\ref{sysV}$), and $W$ be such that $[V,W] =0 $. If $\displaystyle{ \frac 1\xi }$ satisfies the condition $(B)$, then the system $(\ref{sysxi})$ has at most $n$ critical cycles. 
\end{corollary}
{\it Proof.}  
It is an immediate consequence of corollary \ref{cordern}, choosing $\displaystyle{ m = \frac 1\xi }$.
\hfill$\clubsuit$\medbreak

It is not difficult to find examples of systems (\ref{sysxi}) with exactly $n$ critical cycles. It is sufficient to consider a linear center
\begin{equation} \label{isolinear}
x' = y\, \xi(x,y), \qquad y' = -x\, \xi(x,y),
\end{equation} 
with $\xi(z) = \omega(|z|)$,  $\frac 1{\omega(r)}$ having exactly $n$ critical points. In fact, in this case one has $\displaystyle {T(z) = \frac 1{ \xi(z)}} = \frac 1{  \omega(|z|)}$. 

Replacing again the notation $\partial_W \mu $ with $\mu'$ one has
$$
\mu_1 = \mu  =  \frac{ \partial_W m}{m} = \xi \ \partial_W \left(  \frac 1\xi  \right)=  - \frac{\xi'}{\xi}.
$$
$$
\mu_2 = \mu^2 + \mu' = \frac {2 \xi'^2 - \xi'' \xi}{\xi^2} .
$$
$$
\mu_3 = \mu^3 + 3 \mu \mu' + \mu'' =  \frac {-6 \xi'^3 + 6 \xi \xi' \xi ''  - \xi^2 \xi'''  }{\xi^3} .
$$
$$
\mu_4 = \mu^4 + 6 \mu^2 \mu' + 4 \mu \mu'' + 3 \mu'^2  +\mu''' = 
\frac { 24 \xi'^4 - 36 \xi \xi'^2 \xi'' +8 \xi^2 \xi' \xi'''  +6 \xi^2 \xi''^2 - \xi^3 \xi^{(4)}}{\xi^4}.
$$
$$
\mu_5 = \mu^5 + 10 \mu^3 \mu' + 10 \mu^2 \mu'' +15 \mu \mu'^2 +5 \mu \mu''' +10 \mu' \mu'' + \mu^{(4)} = 
$$
$$
= \frac{-120 \xi'^5+240\xi'^3 \xi'' \xi -60 \xi'^2 \xi''' \xi^2 - 90 \xi' \xi''^2 \xi^2 +10 \xi'  \xi^{(4)} \xi^3 +20 \xi''\xi'''\xi^3 -  \xi^{(5)} \xi^4}{\xi^5}.
$$
Let us assume $\xi(x,y)$ to be analytic, $\xi(x,y) = \sum_{n=0}^\infty \xi_n(x,y)$, where $\xi_n =\sum_{j=0}^n c_{n,j} x^{n-j} y^j$ is an $n$-degree homogeneous polynomial. If $\xi > 0$ in a neighbourhood of $O$, then the origin is a center of (\ref{isolinear}). 
One has $V(x,y) =(y,-x)$, $W(x,y) =(x,y)$, so that, by Euler's formula, 
$$
\xi' = \partial_W \xi =  x\xi_x +y\xi_y = x\left( \sum_{n=0}^\infty \xi_n(x,y) \right)_x + y\left( \sum_{n=0}^\infty \xi_n(x,y) \right)_y = \sum_{n=0}^\infty n\, \xi_n(x,y).
$$ 
Similarly
$$
\xi^{(k)} = \partial_W^{(k)}  \xi = \sum_{n=0}^\infty n^k \xi_n(x,y).
$$
Then, working as in corollary \ref{muconv}, if
$$
\xi''  =  \partial_W^{(2)}  \xi\sum_{n=0}^\infty n^2 \xi_n(x,y) \leq 0,
$$
then $T$ is $W$-convex, since $\displaystyle{  \mu_2 = \mu^2 + \mu' = \frac {2 \xi'^2 - \xi'' \xi}{\xi^2} } \geq 0$.
Moreover, if $\xi''$ satisfies the condition $B$, then $T$ is strictly $W$-convex. This is the case of the system
$$
x' = y\, ( 1 - \xi_n(x,y)), \qquad y' = -x\,  ( 1 - \xi_n(x,y)),
$$
if  $\xi_n(x,y) \geq 0$, $\xi_n(x,y)$ homogeneous of degree $n$. On the other hand, since $- \xi' =  n \xi_n(x,y) \geq 0$, $T$ is increasing, hence there are no critical orbits. In this case the period annulus is contained in the oval defined by the inequality $\xi_n(x,y) < 1$. 

In next corollary we consider a class of systems admitting critical orbits.

\begin{corollary}  Let $\xi_l(x,y)$ homogeneous of degree $l$ for $ l=k,n$. Assume $k < n$,  $(3 - 2\sqrt{2}) n \leq k \leq (3 + 2\sqrt{2}) n$, $\xi_k(x,y) >0$ for $(x,y) \neq (0,0)$. Then $O$ is  a center of
$$
x' = y\, ( \xi_k(x,y) - \xi_n(x,y)), \qquad y' = -x\, ( \xi_k(x,y) - \xi_n(x,y)),
$$
 with period function $T$ $W$-convex on $N_O$. If, additionally,  $ \xi_n(x,y) < 0$, then there exists exactly one critical orbit in $N_O$.
\end{corollary}
{\it Proof.}  
The origin is  a center, since the condition $\xi_k >0$ implies that in a neighbourhood of $O$ the orbits coincide with those of the linear center. 
The numerator of 
$$
\mu_2 = \mu^2 + \mu' = \frac {2 \xi'^2 - \xi'' \xi}{\xi^2} 
$$
is
$$
2 \xi'^2 - \xi'' \xi = k^2 \xi_k^2 + (k^2+n^2-4kn) \xi_k \xi_n +n^2 \xi_n^2.
$$
Considering it as a quadratic form in the indeterminates $\xi_k$, $\xi_n$, a sufficient condition for $\mu_2$ not to change sign is
$$
\Delta = (k^2+n^2-4kn)^2 - 4k^2n^2 = (n^2-6nk+k^2)(n-k)^2 \leq 0.
$$
If $(3 - 2\sqrt{2}) n \leq k \leq (3 + 2\sqrt{2}) n$, then $n^2-6nk+k^2 \leq 0$, hence $\Delta \leq 0$.

If,  $ \xi_n(x,y) < 0$, then the central region $N_O$ is contained in the oval having equation  $\xi_k(x,y) - \xi_n(x,y) = 0$. Such an oval consists of critical points. The system is analytical, hence the boundary $\partial N_O$ cannot be a limit cycle, so that at least one of such critical points lies on the external boundary of $N_O$. As a consequence, $T$ goes to $\infty$ approaching the oval. Moreover, $T$ goes to $\infty$ as well approaching $O$, since the center is degenerate. Hence $T$ has a minimum in $A$, reached on a  critical orbit $\phi_V(t,z)$, which is the unique critical orbit in $A$ by $T$'s $W$-convexity.
\hfill$\clubsuit$\medbreak

\subsection{Jacobian maps and Hamiltonian systems with separable variables}  \label{varsep}

Let $\Psi(z) = (P(z),Q(z)) \in C^\infty(\Omega,\R^2)$, with jacobian matrix $J_\Psi (z)$. If  ${\delta}(z) = \det J_\Psi (z) \neq 0$, we say that it is a {\it jacobian map}. If $\Psi$ is a jacobian map, possibly exchanging $P$ and $Q$, we may assume its jacobian determinant ${\delta}(z) = \det J_\Psi (z)$ to be positive on $\Omega$. This ensures local, but not global invertibility. Let us consider the function $H(z)$ defined by  $2H(z) =|\Psi(z)|^2$. 
The Hamiltonian system having $H$ as Hamiltonian is
\begin{equation}\label{syshamPQ}
x' =  H_y =  PP_y + QQ_y, \qquad y' = - H_x =  - PP_x - QQ_x. 
\end{equation} 
We consider also the system obtained dividing (\ref{sysV}) by $\delta(z)$, 
\begin{equation}\label{sysVd}
x' = \frac{ PP_y + QQ_y}{\delta }, \qquad y' = - \frac{ PP_x + QQ_x}{\delta }, 
\end{equation}
and the system
\begin{equation}\label{sysWd}
x' = \frac{ PQ_y - QP_y}{\delta }, \qquad y' = \frac{ -PQ_x + QP_x}{\delta }. 
\end{equation}
From now on, we denote by $V_\Psi$ the vector field of (\ref{sysV}), $V_\delta  $ the vector field of (\ref{sysVd}), $W_\delta  $  the vector field of (\ref{sysWd}).  We denote by $\phi_{V_\Psi}(t,z)$ a solution to (\ref{sysV}), and by $\phi_{W_ \delta}(t,z)$ a solution to  (\ref{sysWd}). 

In \cite{Vo} critical points of (\ref{sysVd}) were proved to be isochronous centers, under the assumption $\delta > 0$. In fact, the systems 
(\ref{sysVd})  and (\ref{sysWd}) commute with each other, that allows us to apply the theorem \ref{teordern} to find a formula for the $n$-th derivative of some hamiltonian period functions. We denote by $\partial_H^{(n)} T$ the $n$-th derivative of $T$ with respect to the function $H$.

\begin{theorem} \label{teorPQ} Let $\Psi(z) = (P(z),Q(z)) \in C^\infty(\Omega,\R^2)$.  If $A$ is a period annulus of (\ref{syshamPQ}), then 
\begin{equation}  \label{derjac}
\partial_{W_\delta}^{(n)} T(z) = \int_0^{T(z)} \delta(\phi_{V_\Psi}(t,z))  \  \partial_{W_\delta}^{(n)}\left(  \frac 1{\delta(\phi_{V_\Psi}(t,z)) } \right) \, dt = 
\end{equation}
\begin{equation}  \label{derjacmu}
= \int_0^{T(z)}  \mu_{\psi n}(\phi_{V_\Psi}(t,z))\, dt,
\end{equation}
where $\mu_{\psi n}$ is defined recursively as in theorem \ref{ricorr}, with $\displaystyle{\mu_{\psi} = -  \ln \delta }$.
Moreover,
\begin{equation}  \label{derjacPQH}
\partial_{H}^{(n)} T(z) = \int_0^{T(z)}  \overline \mu_n(\phi_{V_\Psi}(t,z))\, dt,
\end{equation}
where $\overline \mu_n$ is defined recursively as in theorem \ref{ricorr}, with $\displaystyle{\overline \mu = - \frac { \ln \delta  }{2H} }$.
\end{theorem}
{\it Proof.} 
The transformation $\Psi$ takes locally the systems (\ref{sysVd}) and (\ref{sysWd}) into the systems (\ref{sysuv}), which commute with each other. Since the system (\ref{syshamPQ}) is obtained from (\ref{sysVd}) multiplying by $\delta$, setting $ \displaystyle{m= \frac 1 \delta} $ one has $[m V_\Psi,W_\delta] = 0$. 
Then the equality (\ref{derjac}) comes from theorem \ref{teordern}.
By lemma \ref{lemmam}, one has
$$
\mu_\Psi  = \frac{\partial_{W_\delta} m}{m} = \delta\  \partial_{W_\delta}\left(  \frac 1\delta \right) =  \delta\  \left(-  \frac {\partial_{W_\delta} \delta}{\delta^2 }\right) =  - \partial_{W_\delta } (\ln\delta ).
$$
The formula (\ref{derjacmu}) comes from theorem  \ref{ricorr}.

The function $H$ vanishes only at critical points of $V_\Psi$, which coincide with those of $V_\delta$ and $W_\delta$, hence the vector field $\displaystyle{  \frac { W_\delta }{2H} }$ is a non-trivial normalizer of (\ref{syshamPQ}), with N-cofactor  $\displaystyle{\overline \mu = - \frac { \ln \delta  }{2H} }$. The formula (\ref{derjacPQH}) comes as well from theorem \ref{ricorr}.
\hfill$\clubsuit$\medbreak
 
In theorem's \ref{teorPQ} proof one does not need the invertibility of $\Psi$ on all of the period annulus, since one only needs a C-factor, provided by $\displaystyle{ \frac 1{\delta}  }$. 

The applicability of theorem \ref{teorPQ} depends on the possibility to write a  hamiltonian function $H(z)$ in the form $2H(z) =|\Psi(z)|^2$, with $\det J_\Psi(z) \neq 0$.  Such a question was addressed in \cite{SM,MRT}.
 
Since the above theorem can be stated in terms of normalizers, we can consider also non-hamiltonian systems, as observed in section 2. 
 
 \begin{corollary} \label{corPQ} 
 Let $\rho > 0$ be an inverse integrating factor of the system (\ref{sysV}), with first integral $H(x,y)$ satisfying $2H(z) =|\Psi(z)|^2$, for some jacobian map $\Psi$. If $A$ is a period annulus of (\ref{syshamPQ}), then 
\begin{equation}  \label{dercorPQ}
\partial_{W_\delta}^{(n)} T(z) = \int_0^{T(z)}  \mu_n(\phi_{V_\Psi}(t,z))\, dt,
\end{equation}
where $\mu_n$ is defined recursively as in theorem \ref{ricorr}, with $\mu = -  \ln (\delta \rho )$.
Moreover,
\begin{equation}  \label{dercorPQH}
\partial_{H}^{(n)} T(z) = \int_0^{T(z)}  \overline \mu_n(\phi_{V_\Psi}(t,z))\, dt,
\end{equation}
where $\overline \mu_n$ is defined recursively as in theorem \ref{ricorr}, with $\displaystyle{ \overline \mu = - \frac { \ln (\delta \rho ) }{2H} }$.
\end{corollary}
{\it Proof.} 
By hypothesis, one has $(-V_2(z),V_1(z) )= \rho(z) \nabla H(z)$, with $2H(z) =|\Psi(z)|^2$, hence $V(z) = \rho(z) V_\Psi(z)$. Working as in section 2, one has
$$
[V, {W_\delta } ] =[\rho V_\Psi, {W_\delta } ] =  ( \mu_\Psi \rho - \partial_{W_\delta } \rho) V =
 \Big( - \partial_{W_\delta } (\ln\delta ) - \partial_{W_\delta }  ( \ln \rho) \Big) (\rho V) =
$$ $$
  \Big( - \partial_{W_\delta } \big(\ln (\delta \rho) \big)  \Big) (\rho V_\Psi) =  \Big( - \partial_{W_\delta } \big(\ln (\delta \rho) \big)  \Big) V .
$$
that gives the formula (\ref{dercorPQ}).  The formula (\ref{dercorPQH}) can be proved as formula (\ref{derjacPQH}) in theorem \ref{teorPQ}.
 \hfill$\clubsuit$\medbreak

In this paper we only consider the applicability of theorem \ref{teorPQ}  to some hamiltonian systems with separable variables. Let $I_F$, $I_G$ be intervals containing $0$, $F\in C^\infty(I_F,\R)$, $ G \in C^\infty(I_G,\R)$. Assume $G(x)$ and $F(y)$ to have isolated minima at $0$. Then the origin $O$ is a center of the  hamiltonian  system having $H(x,y) = G(x) + F(y)$ as hamiltonian function,
\begin{equation} \label{varsep}
\left\{  \begin{array}{rl}
x' &= F'(y) \\
y' &= - G'(x) .
\end{array}   \right.   
\end{equation} 
Under quite general conditions, such systems can be considered as special cases of (\ref{syshamPQ}), obtained taking $\Psi(x,y) = (P(x), Q(y) ) = (s(x)\sqrt{2G(x)}, s(y)\sqrt{2F(y)})$, where $s(t)$ the sign function, assuming values $-1,0,1$ for $t<0$, $t=0$, $t>0$, respectively. The jacobian determinant of such a map is $\delta(x,y) =  P'(x) Q'(y) =$
$\displaystyle{ s(xy)  \frac{F'(y)G'(x)}{2\sqrt{F(y)G(x)}}}$. In this case, the normalizer (\ref{sysWd}) has the form
\begin{equation}\label{sysWdvarsep}
x' = \frac{ PQ_y - QP_y}{\delta } =  \frac{ P}{P' }=  \frac{ 2G}{G' }, \qquad y' = \frac{ -PQ_x + QP_x}{\delta } =  \frac{ Q}{Q' }=  \frac{ 2F}{F' }, 
\end{equation}
and differs from that one given in \cite{FGG} for the presence of a factor 2. Also the corresponding  N-cofactor 
$$
\mu_{s}(x,y) = - \partial_{W_\delta } (\ln \delta(x,y) ) = - \partial_{W_\delta } \Big(\ln (P'(x) Q'(y)  \Big) = 
$$
$$
= - \partial_{W_\delta } \left(\ln  \frac{s(x)G'(x)}{\sqrt{2G(x)}} \right)  - \partial_{W_\delta } \left(\ln  \frac{s(y)F'(y)}{\sqrt{2F(y)}} \right) =
  2\left(  \frac{G(x)}{G'(x)} \right)' + 2\left( \frac{F(y)}{F'(y)} \right)'   -2  
$$
differs from that one given in \cite{FGG} for the presence of a factor 2.
The commuting system (\ref{sysVd}) has the form
\begin{equation} \label{isovarsep}
x' =  \frac{2s(xy) \sqrt{G(x)F(y)}}{G'(x)} = \frac{Q(y)}{P'(x)}, \qquad y' = -  \frac{ 2s(xy)\sqrt{G(x)F(y)}}{F'(y)} = -\frac{P(x)}{Q'(y)}.
\end{equation}  
Let us set 
$$
m_s(x,y) = \frac{1}{\delta(x,y)} = \frac{2s(xy)\sqrt{G(x)F(y)}}{G'(x)F'(y)},
$$
whenever $F(y)G(x) \geq 0$, $F'(y)G'(x) \neq 0$. We  say that a function $L:I_L\rightarrow \R$, $I_L$ open interval containing 0, satisfies the hypothesis $(Sep)$ if 
\begin{itemize}
\item[$sep_1)$] $L\in C^\infty(I_L,\R)$, $L(0) =0$, $tL'(t) > 0$ in $I_L$, $L$ is not flat at 0;
\item[$sep_2)$] $\displaystyle{\frac{\sqrt{L}}{L'}\in C^\infty(I_L,\R)}$, $s(t)\sqrt{L(t)} \in C^\infty(I_L,\R)$. 
\end{itemize}
The condition $sep_1)$ implies $L(t) > 0$ for $t$ close to 0, $t \neq 0$, hence its Taylor expansion starts with an even power of $t$,
$$
L(t) = c_{2k} t^{2k} + o(t^{2k}).
$$ 
As a consequence,
$$
L'(t) = 2kc_{2k} t^{2k-1} + o(t^{2k-1}),
$$ 
so that
$$
\frac{\sqrt{L}}{L'} = \sqrt{ \frac{L}{L'^2} } =  \sqrt{ \frac{c_{2k} t^{2k} + o(t^{2k})}{{4k^2c_{2k}^2 t^{4k-2} + o(t^{4k-2})}} } .
$$ 
Such a function can be differentiable only if $4k-2 \leq 2k$, i. e. $k \leq 1$, hence $k=1$. This agrees with the fact that non-degeneracy is a necessary condition for a center's linearizability \cite{V}. 

\begin{theorem}  \label{teorvarsep} Let $F$ and $G$ satisfy $(Sep)$ on $I_F$, $I_G$, open intervals containing 0. Then
\begin{itemize}
\item[i)] the transformation $\Psi(x,y) = \left(  s(x)\sqrt{2G(x)}, s(y)\sqrt{2F(y)} \right)$ linearizes both systems (\ref{sysWdvarsep}) and (\ref{isovarsep});
\item[ii)] $O$ is an isochronous center of the system (\ref{isovarsep});
\item[iii)] the period function of the system (\ref{varsep}) satisfies
\begin{equation} \label{dernvarsep}
\partial_{W_\delta}^{(n)} T(z) = \int_0^{T(z)} \frac {\partial_{W_\delta}^{(n)} m_s(\phi_{V_\Psi}(t,z)) }{m_s(\phi_{V_\Psi}(t,z))} dt,
\end{equation}
\end{itemize}
where the integration is performed along the cycle $\phi_{V_\Psi}(t,z)$ starting at $z$.
\end{theorem}
{\it Proof.} One can prove $i)$ by direct computation. The systems (\ref{isovarsep}) and (\ref{sysWdvarsep}) are transformed into the systems (\ref{sysuv}). 

In order to prove $ii)$, it is sufficient to observe that the systems (\ref{sysuv}) commute, hence also (\ref{isovarsep}) and (\ref{sysWdvarsep})  commute. By the main result in \cite{S0}, every period annulus of the system (\ref{isovarsep}) is isochronous. 

The statement $iii)$ is an immediate consequence of theorem \ref{teorPQ}. 
\hfill$\clubsuit$\medbreak

The theorem \ref{teorvarsep} does not apply to degenerate centers of hamiltonian systems. In order to study such systems' higher order derivatives, one can apply the corollary  \ref{ricorr}, starting with a suitable normalizer.  

A "universal" normalizer, as pointed out in section 2 (see \cite{S4}), that allows to cover both centers and period annuli, is 
$$
x' = \frac{H_x}{H_x^2 + H_y^2} =  \frac{G'}{G'^2+F'^2}, \qquad y' = \frac{H_y}{H_x^2 + H_y^2} =  \frac{F'}{G'^2+F'^2},
$$
which gives the N-cofactor
\begin{equation}\label{muhamsep}
\mu_H =   \frac{( F'' - G'')(G'^2  - F'^2)}{(G'^2+F'^2)^2} .
\end{equation}
This allows to cover both degenerate and non-degenerate hamiltonians, obtaining  $T$'s derivatives without need to explicitly  find the corresponding $m$, by using the recursive formulae of theorem  \ref{ricorr}. Moreover, such derivatives are computed with respect to the variable $H$.  On the other hand, even in the simplest non-degenerate cases, this leads to more complex computations than using $m_s$ and the related normalizer $\mu_{s}$. As an example, If $\displaystyle { H(x,y) = \frac{x^2 + x^4 + y^2}2 }$, one has
$$
\mu_{s} = - \frac{x^2 (3+2x^2)}{(1+2x^2)^2}
$$ 
$$
\mu_H = -  \frac{6x^2(x^2+4x^4+4x^6-y^2) }{(x^2+4x^4+4x^6+y^2)^2} .
$$
\medbreak
The same N-cofactor can be obtained by using the system
$$
x' = \frac{G}{H G'}, \qquad y' = \frac{F}{H F'},
$$
which is defined on all of the central region $N_O$, but is not defined in a neighbourhood of every point $(x,y)$ where $G'(x) = 0$ and $G(x) \neq 0$, or  $F'(y) = 0$ and $F(y) \neq 0$.  This is the case of the hamiltonian system
$$
x' = y, \qquad y' = -x+2x^2-x^3.
$$
One has
$$
\mu_{s} = \frac{x (3x^2-9x+8) }{12(1-x)^3},
$$
which diverges at 1. The hamiltonian system  has a center at $O$, with central region defined by 
$\displaystyle{ N_O =\left\{ (x,y) :0 < \frac {x^2}2  - \frac {2x^3} 3 + \frac {x^4} 4 + \frac {y^2} 2  < \frac 1{12} \right\} 
}$. 
The central region's external boundary consist of a critical point at $(1,0)$ and a homoclinic orbit having $(1,0)$ as $\alpha$- and $\omega$-limit set. Every orbit external to $N_O$ is a cycle enclosing $N_O$. One cannot study $T$'s derivatives on such cycles by means of $\mu_{s}$, since every external cycle has a point (actually, two points) where $\mu_s$ is not defined.
Obviously, also the jacobian map approach can no longer be applied, since $\delta$ vanishes where $G'$ vanishes. On the other hand, we can compute $\mu_H$ for the normalizer
$$
x' = \frac{x-2x^2+x^3}{(x-2x^2+x^3)^2+y^2},  \qquad y' =\frac{y}{(x-2x^2+x^3)^2+y^2} .
$$
One has
$$
\mu_H = \frac{x(-4+3x)(y+x-2x^2+x^3)(-y+x-2x^2+x^3)}{(x-2x^2+x^3)^2+y^2}.
$$
The study if $T$'s monotonicity requires an additional study of $\int_0^T \mu_H$, in order to determine its sign. Higher order derivatives require the study of much more complex functions.

\medbreak

When dealing with the period function of an analytic system (\ref{varsep}) in a central region $N_O$, the second, and more convenient approach, consists in applying the theorem \ref{ricorr} to $\mu_s$. Assume 
$$
G(x) = O( x^{2k} ), \qquad F(y) = O(y^{2h}). 
$$
Then the system (\ref{sysWdvarsep}) is of class $C^\infty$ on $N_O$, as well as $\mu_s$. Applying the  theorem \ref{ricorr} does not require to produce the  corresponding C-factor, even if in some simple cases it can be found, as for the systems
$$
x' = y^{2k-1}, \qquad y' = - x^{2k-1}, \qquad k \in \N, \quad k >0.
$$
The normalizer (\ref{sysWdvarsep}) is
$$
x' = \frac xk, \qquad y' = \frac yk,
$$
with N-cofactor $\displaystyle{  \mu_s = \frac {2(k-1)}{k}  }$ and $\displaystyle{  \mu_{s2} = \mu_s^2 = \frac {4(k-1)^2}{k^2} }$. There exist infinitely many C-factors, all producing $\mu_s$ as N-cofactor. Among them, one has all the functions  $\displaystyle { m(x,y) = \frac {1}{\nu(x,y)}}$, with $\nu$ homogeneous of degree $2k-2$.

\medbreak
In order to perform next computations with the N-cofactor $\mu_{s}$, we find convenient to write next formulae in terms of the functions $\displaystyle{\Gamma(x) = \frac{2G(x)}{G'(x)}}$, $\displaystyle {\Phi(y) =  \frac{2F(y)}{F'(y)} }$, using the normalizer (\ref{sysWdvarsep})
$$
x' = \Gamma(x), \qquad y' = \Phi(y),
$$
as well as its N-cofactor
$$
\mu_{s} (x,y) =  \Gamma'(x) +  \Phi'(y)   - 2.
$$
One can write $\mu_{s} $ as follows, choosing $\alpha$ and $\beta$ such that $\alpha + \beta = 1$.
$$
\mu_{s}  = \left( \Gamma'(x)  -\alpha \right) +  \left( \Phi'(y) - \beta  \right) = 
\frac{(2-\alpha )G'^2 - 2GG''}{G'^2}  + \frac{(2-\beta) F'^2 - 2FF''}{F'^2}.
$$
Hence, if there exists a couple $(\alpha,\beta)$ such that on a period annulus $(2-\alpha )G'^2 - 2GG'' \geq 0$, and $(2-\beta) F'^2 - 2FF'' \geq 0$, then $T$ is increasing. For systems equivalent to second order conservative equations one does not have such a freedom of choice, since 
$$
F(y) = \frac{y^2}2, \quad \Phi'(y) = 1, \quad  \mu_{s}  =  \left( \frac{2G(x)}{G'(x)} \right)'   - 1 = \frac{G'^2-2GG''}{G'^2}.
$$
As a consequence, if $G'^2-2GG'' \geq 0$ in a period annulus $A$, then $T$ is increasing in $A$ (see \cite{FGG}). 

\medbreak

Let us denote by $W_s$ the  normalizer vector field (\ref{sysWdvarsep}). One can compute $\mu_{s2} = \mu_s^2 + \partial_{W_s} \mu_s$:
$$
\mu_{s2} = 4 \bigg[  1   +  2\,\frac{ G  G'' }{G'^2 } \cdot\frac{ F  F''}{F'^2 } +
$$ $$
+  \frac { 3 G^2 G''^2 - 3 G G'^2 G''  -  G^2 G' G''' }{G'^4} + \frac{ 3  F^2 F''^2 - 3 F F'^2 F'' -   F^2 F' F''' }{F'^4} \bigg].
$$
Such a formula has been written in such a way to emphasize the presence of a single term with mixed variables.

If $O$ is a center of (\ref{varsep}), we denote by   $N_O^x$,  $N_O^y$ the projections of $N_O$ on the $x$-axis and on the $y$-axis, respectively.

\begin{corollary} \label{mu2varsep} Assume $G(x) =b_{2k} x^{2k}  +  r_g(x) \in C^\infty(I_G,\R) $, $ F(y) =  c_{2h} y^{2h}  +  r_f(y) \in C^\infty(I_F,\R)$, $0 \in I_G \cap I_F$, $0 <k,h \in \N$, $b_{2k}, c_{2h} > 0$, $r_g(x) = o(x^{2k})$, $r_f(y) = o(y^{2h})$. Then $O$ is a center, (\ref{sysWdvarsep}) is of class $C^ \infty(N_O,\R^2)$,  $\mu_s\in C^ \infty(N_O,\R^2)$.
 If $ \mu_{s2} \geq 0$ ($\leq 0$) in $N_O^x \times N_O^y$, then $T$ is $W_\delta$-convex (concave) on $N_O$. If, additionally,  on every cycle of $N_O$ there exists a point where such  inequality holds strictly, then $T$ is strictly $W_\delta$-convex (concave).
\end{corollary}
{\it Proof.} 
The functions $F'(y)$ and $ G'(x)$ vanish only at $0$. As for the regularity of (\ref{sysWdvarsep}) and $\mu_s$ at $x=0$ and $y=0$, one has
$$
\frac{G(x)}{G'(x)} = x \ \frac{ b_{2k} +   { \overline r}_g(x) }{ 2k\, b_{2k} +   \overline { r'}_g(x)  }, 
\quad  { \overline r}_g(x) = \frac {r_g(x)}{x^{2k}}, \quad  \overline { r'}_g(x) = \frac {r'_g(x)}{x^{2k-1}}
$$
which is the product of two functions of class $C^\infty$ at $O$. Similarly for $F(y)$. Since $\mu_n$ is obtained multiplying $\mu$'s derivatives  by the components of (\ref{sysWdvarsep}), one has $\mu_n \in C^\infty$.
Then the (strict) $W_\delta$-convexity comes from the formula 
$$
\partial_{W_\delta}^{(2)} T(z) = \int_0^{T(z)} \mu_{s2}(\phi_{V_\Psi}(t,z)) \ dt.
$$
\hfill$\clubsuit$\medbreak

Looking for simpler conditions which imply $T$'s convexity, one can write
$$
\mu_{s2} =  ( \Gamma' + \Phi' - 2 )^2 + \Gamma'' \Gamma +  \Phi'' \Phi = $$ $$
= ( \Gamma' + \Phi')^2 - 4 \Gamma' - 4 \Phi' +  4  + \Gamma'' \Gamma +  \Phi'' \Phi.
$$
Here $\mu_{s2}$ appears in two ways as the sum of a squared term plus terms depending only on $x$ or $y$. This allows to study its sign in a simpler way than with genuine two-variables functions. We say that a function $L(t)$ satisfies the condition $(C_\lambda)$ if
$$
\zeta_ \lambda = 2 L^2 L''^2 +L L'^2 L''  - L^2 L' L''' \geq  \lambda L'^4, \qquad \lambda \in \R.
$$
We say that a function $L(t)$ satisfies the condition $(C)$ if
$$
\zeta =  2 L L''^2 - L'^2 L''   - L L' L''' \geq 0.  
$$

\begin{corollary} \label{mu2varsepridotto} 
Under the hypotheses of corollary  \ref{mu2varsep}, if one of the following condition holds, then $T$ is $W_s$-convex.
\begin{itemize}
\item[i)] $G$ satisfies $(C_\alpha)$ on   $N_O^x$, $F$ satisfies $(C_\beta)$  $N_O^y$, respectively, with $\alpha + \beta - 3 \geq 0$;
\item[ii)] both $G$ and $F$ satisfy the condition $(C)$ on   $N_O^x$ and  $N_O^y$, respectively;
\end{itemize}
\end{corollary}
{\it Proof.} 
i) One has
$$
\mu_{s2i} =  - 4 \Gamma' - 4 \Phi' +  4 + \Gamma'' \Gamma +  \Phi'' \Phi = 
$$ 
$$
4 \bigg[  \frac { 2 G^2 G''^2 +G G'^2 G''  - G^2 G' G'''}{G'^4} + 
\frac { 2 F^2 F''^2 + F F'^2 F'' - F^2 F' F''' }{F'^4} - 3 \bigg] \geq 
$$
$$
\geq 4 (  \alpha + \beta  - 3) \geq 0.
$$
hence $\mu_{s2} = ( \Gamma' + \Phi')^2 - 4 \Gamma' - 4 \Phi' +  4  + \Gamma'' \Gamma +  \Phi'' \Phi \geq 0$.

ii) One has
$$
\mu_{s2ii} =\Gamma'' \Gamma +  \Phi'' \Phi = 4 \bigg[  \frac {2 G^2 G''^2 - G G'^2 G''   - G^2 G' G'''}{G'^4} + 
\frac {2 F^2 F''^2 - F F'^2  F'' -  F^2 F' F'''}{F'^4}  \bigg].
$$
Since 
$$
  \frac {2 G^2 G''^2 - G G'^2 G''   - G^2 G' G'''}{G'^4}  =   \frac {G}{G'^4} \bigg(  2 G G''^2 - G'^2 G''   - G G' G''' \bigg),
$$
and $   G  \geq 0 $ in $N_O$, under the the condition $(C)$ one has $\Gamma'' \Gamma \geq 0$ in $N_O$. The same holds for $\Phi'' \Phi $, hence $\mu_{s2} =  ( \Gamma' + \Phi' - 2 )^2 + \Gamma'' \Gamma +  \Phi'' \Phi \geq 0$. 
\hfill$\clubsuit$\medbreak

An important special case is that of second order conservative ODE's. In this case  $F(y) = \displaystyle{ \frac {y^2}2   }$. In next corollary we just write the conditions $i)$ and $ii)$ of corollary \ref{mu2varsepridotto} specialized to such a case.

\begin{corollary} \label{mu2ODE} 
Under the hypotheses of corollary  \ref{mu2varsep}, if $F(y) = \frac {y^2}2   $, one has 
$$
\mu_{s2} = \frac{G'^4 - 8 G'^2 G G'' + 12 G^2 G''^2 -4 G^2 G''' G'}{G'^4} .
$$
Moreover, the condition $i)$ of corollary  \ref{mu2varsepridotto} holds if and only if 
$$
G G'^2  G'' + 2 G^2 G''^2 - G^2  G' G''' - 2 G'^4 \geq 0.
$$
The condition  $ii)$ of corollary  \ref{mu2varsepridotto} does not change.
\end{corollary}

We can apply the corollary  \ref{mu2varsepridotto} ii) to the system
\begin{equation}  \label{sys1}
x' = 2 \sin y \cos y, \qquad y' = -2 \sin x \cos x,
\end{equation}
obtained taking $G(x) = (\sin x)^2$, $F(y) = (\sin y)^2$. The origin is a center whose central region is a punctured open square having vertices at $\left(\frac \pi 2, 0\right)$, $\left(0, \frac \pi 2\right)$, $\left(-\frac \pi 2, 0\right)$, $\left(0,-\frac \pi 2\right)$. One has $\Gamma(x) = \tan x$, $\Gamma''(x) \Gamma(x) = 2( \tan^2 x) (1 + \tan^2 x) \geq 0$ in $N_O$. Similarly for $\Phi(y)$.
\begin{figure}[h!]
  \caption{The central region of system (\ref{sys1}).}
  \centering
\includegraphics[width=5.4cm]{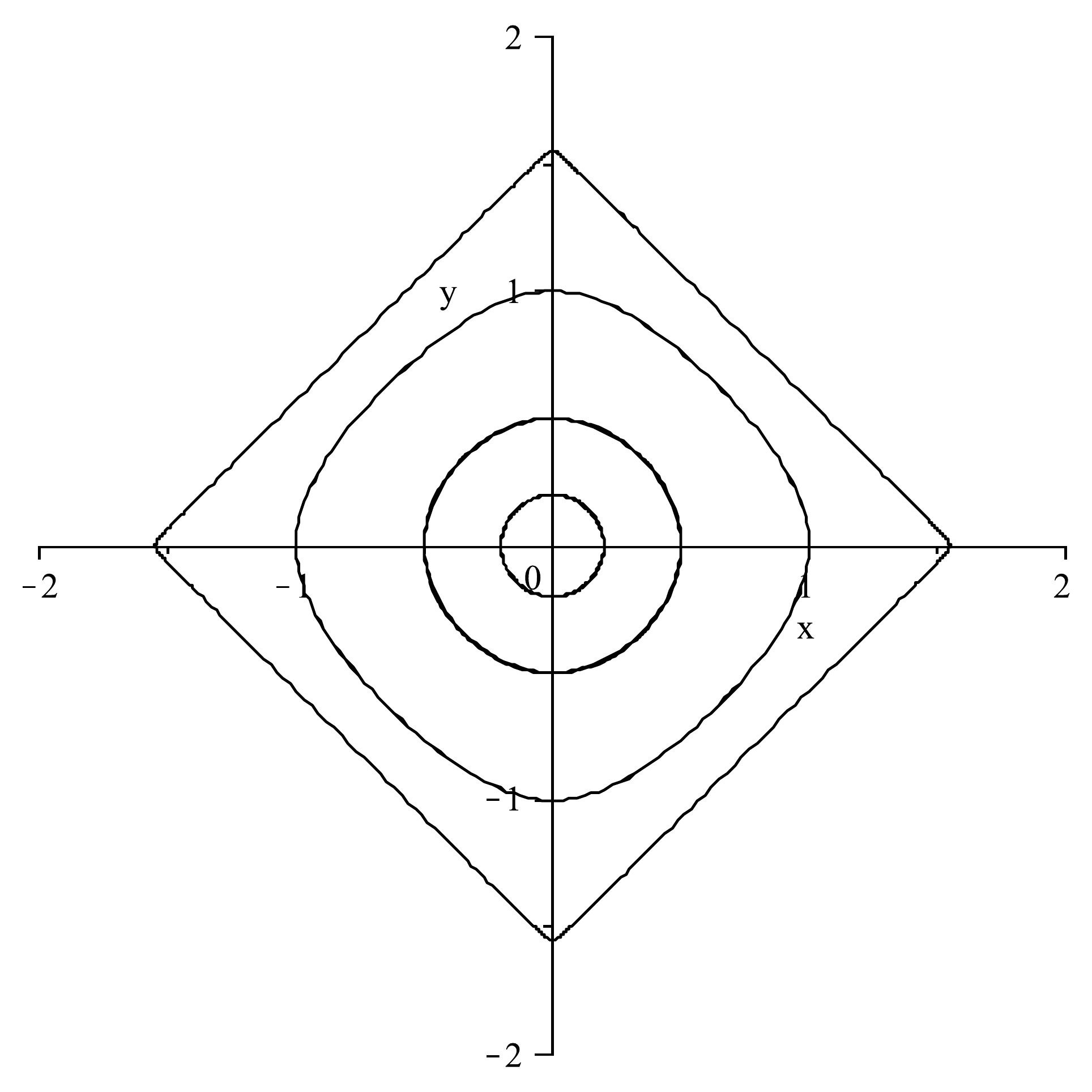}  
\end{figure}

We have an example of system satisfying the hypotheses of corollary  \ref{mu2ODE}  taking 
$\displaystyle{  G(x) = e^{(\cos^3 x - 3 \cos x)}   }$, $\displaystyle{  F(y) = \frac {y^2} 2   }$. The related system is 
\begin{equation}  \label{sys2}
x' = y, \qquad y' = - 3 \sin^3 x e^{(\cos^3 x - 3 \cos x)}.
\end{equation}
The origin is a degenerate center, with bounded central region $N_O$. Its boundary $\partial N_O$ contains a critical point, hence $T$ diverges as a cycle approaches $\partial N_O$ (see Figure 1).  One has $\displaystyle{  \Gamma(x) =\frac 1{3 \sin^3 x}  }$, $\displaystyle{  \Gamma''(x) \Gamma(x) = \frac{1 + 3 \cos^2 x }{3 \sin^8 x} > 0  }$, $\displaystyle{  \Phi''(y) \Phi(y) =  0  }$ in $N_O$, hence $T$ is strictly convex. This implies the existence of a single critical point, since $T$ goes to $\infty$ both as a cycle approaches $O$ and as it approaches  $\partial N_O$.
\begin{figure}[h!]
  \caption{The central region of system (\ref{sys2}).}
  \centering
\includegraphics[width=5.4cm]{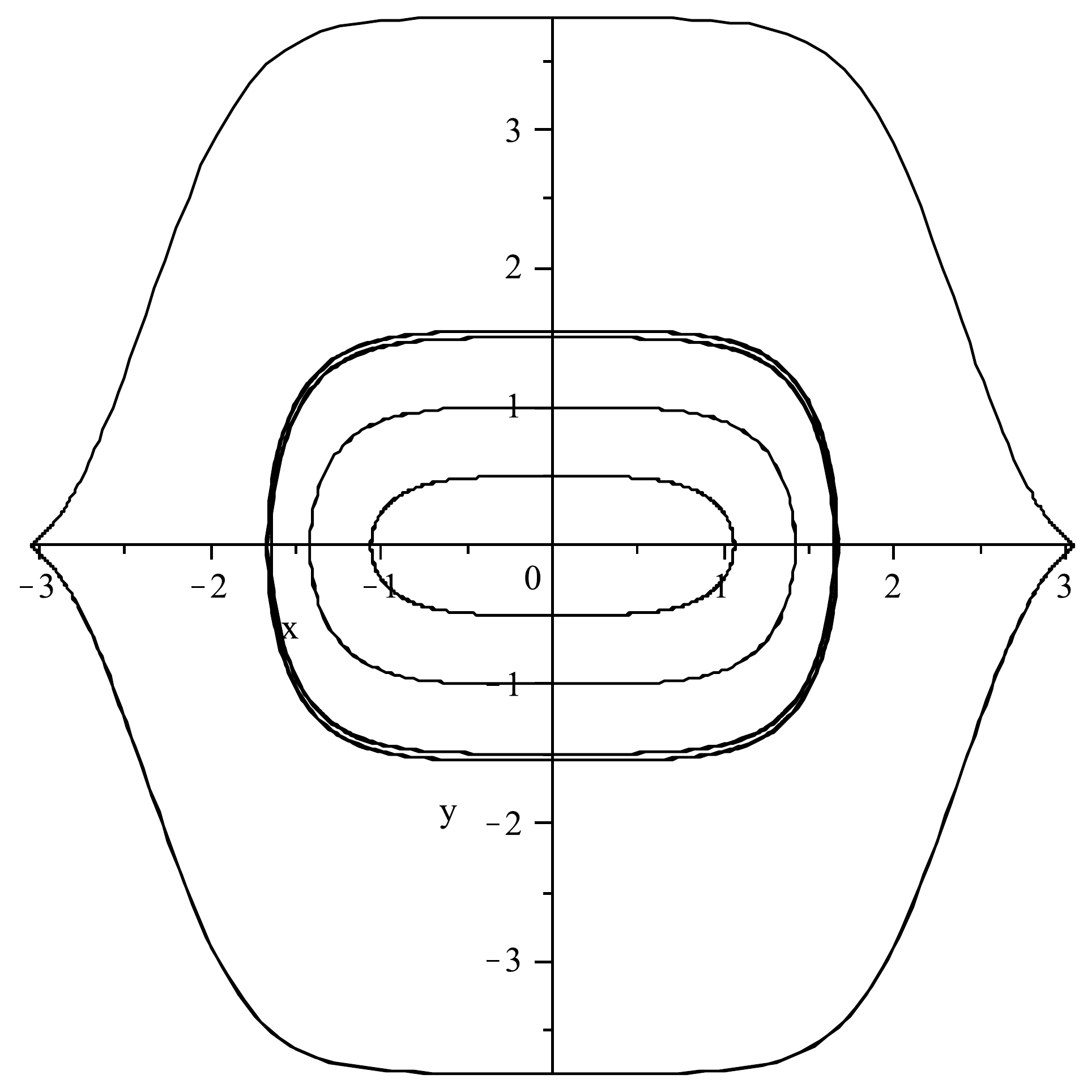}  
\end{figure}

The possibility to consider separately the terms depending only on $x$ and those depending only on $y$ implies that we may swap the components of the above systems, obtaining a new system with convex period function,
$$
x' = 2 \sin y \cos y, \qquad y' = - 3 \sin^3 x e^{(\cos^3 x - 3 \cos x)}.
$$

In \cite{MV,S3}, potential functions of the form $\displaystyle{ \frac {x^m}{ax^n+b}   }$ were considered. We consider here some other rational functions.

\begin{corollary}  \label{potfrac1}
Let $\displaystyle{ G(x) =  \frac {P(x)}{1+P(x)} }$, with $P(x) = ax^4+bx^6+cx^8$. If $a \geq 0$, $b \geq 0$, $c \geq 0$, $abc \neq 0$, 
$20 a^4 c \leq 267 b^2 a c+194 a^2 c^2+18 b^4+90 a^6+219 a^3 b^2$, $8a^3c \leq  40a^5+99a^2b^2$, then $O$ is a center, and $T$ has exactly one critical orbit in $N_O$.
\end{corollary}
{\it Proof.} One has $\displaystyle{ G'(x) =  \frac {P'(x)}{(1 + P(x))^2} } $.
In a neighbourhood of $O$, $xG'(x) \geq 0$,  $G'(x) = 0$ if and only if $x = 0$, hence  $O$ is  a center. 
If $G'(x) = 0$ only at $x=0$, the central region $N_O$ is a strip defined by $y^2 < 2$, since $\displaystyle { \lim_{x \rightarrow \pm \infty} } G(x) = 1$. The vector field is bounded on such a strip, hence $T$ goes to $+\infty$ as the cycle approaches the boundary $\partial N_O$. The same occurs as the cycle approaches the origin, since $O$ is a degenerate critical point.  

If there exists a point $\overline x \neq 0$ such that $G'(\overline x) =0$, then the central region has at least a critical point on its boundary, hence $T$ goes to infinity as the cycle approaches $\partial N_O$. 
Moreover, one has 
$$
\mu =\frac{1}{D(x)}  \bigg(
 36 {x}^{16}{c}^{3}+89 {x}^{14}b{c}^{2}+ \left( 66 a{c}^{2}+74 {b}^
{2}c \right) {x}^{12}+ \left( 110 abc+21 {b}^{3} \right) {x}^{10}
$$ 
$$
+ \left( 48 a{b}^{2}  -12 {c}^{2}+40 {a}^{2}c \right) {x}^{8}
 \left( 37 {a}^{2}b-19 bc \right) {x}^{6}+
$$
$$
+ \left( 10 {a}^{3}-6 {b}^{2}-18 ac \right) {x}^{4}-9 ab{x}^{2}-2 {a}^{2}  \bigg),
$$
with $D(x) = (2a+3bx^2+4cx^4)^2$.  Finally,
$$
\mu_{s2} = \frac {A(x)+B(x)}{D(x)^2}
$$
with $A(x)= (534 b^2 a c+388 a^2 c^2-40 a^4 c+36 b^4+180 a^6+438 a^3 b^2) x^8 + (-8a^3c+40a^5+99a^2b^2)x^4$, $B(x)$ polynomial of degree 32 with positive coefficients, such that $B(0) = 4 a^4 > 0$.
Under the condition given in the hypothesis, $A$ is non-negative, hence $T$ is strictly convex.
\hfill$\clubsuit$\medbreak

Let su consider the following functions:
\begin{itemize}
\item[j)] the functions $\displaystyle  {  \frac{1}{a+ bt^4 + ct^8}  }$, with $a, b, c > 0$, $\displaystyle  {  5b^2 \leq 56 ac  }$;
\item[jj)] the functions $\displaystyle  {  \frac{at^{4k}}{1+bt^{2k}+ct^{6k}}  }$, with $k$ positive integer, $a, b, c > 0$, 
$2  b^3 k \leq + 108 a^2 c k +18 a^2 c + b^3 $.
\end{itemize}

\begin{corollary}  \label{potfrac2}
Let $ G(x) $ and $F(y)$ be of type j) or jj). Then the period function of the corresponding system (\ref{varsep}) is convex.
\end{corollary}
{\it Proof.}  Under the given hypotheses, all the functions of the points j), jj), jjj) satisfy the condition (C). In fact, computing the expression $\zeta = 2 L L''^2 - L'^2 L''   - L L' L'''$ for such functions gives the following.

j) $\displaystyle{  \zeta =  \frac{64 x^4 (6 b c^2 x^{12}+(56 a c^2 -5 b^2 c) x^8+18 a b c x^4 +3 a b^2)}{(a+b x^4+c x^8)^6} \geq 0  }$, if   $a, b, c > 0$, $\displaystyle  {  5b^2 \leq 56 ac  }$;

jj) 
$$ 
\zeta = \frac {A(x)+B(x)}{D(x)^2},
$$
with $D(x) = x^2 ({a+ bx^4 + cx^8})^7$, $B(x)$ polynomial of degree $32 k - 2$ with positive coefficients, $A(x) = (-2 a b^3 k + 108 a^3 c k +18 a^3 c + ab^3 ) x^{18k-2}  + (  54 k - 9) ac^3 x^{30k-2} + ( 16 k - 4) bc^3x^{32k-2}$. Since $k$ is a positive integer, one has $ 54 k - 9 > 0$,  $16 k - 4 >0$. The same holds for the coefficient of $x^{18k-2}$, under the hypothesis $2  b^3 k \leq + 108 a^2 c k +18 a^2 c + b^3 $. 

By corollary  \ref{mu2varsepridotto}, ii), if both $G$ and $F$ satisfy one of the above, one has $T$'s convexity.
\hfill$\clubsuit$\medbreak

We cannot include in the above list the functions $\displaystyle{   \frac {P(x)}{1+P(x)}  }$ of corollary  \ref{potfrac1}, since in that case one has 
$$
2 G^2 G''^2 - G G'^2 G''   - G^2 G' G''' = \frac {-6a^2 b + o(1)}{(1 + ax^4+bx^6+cx^8)^6}. 
$$

\enddocument